\documentclass[12pt,a4paper]{article}
\usepackage[english]{babel}
\usepackage[utf8x]{inputenc}
\usepackage{amsmath}
\usepackage[margin=1in]{geometry}
\usepackage{graphicx,cite}
\usepackage[section]{placeins}
\usepackage{float}
\usepackage{ amssymb }
\usepackage[svgnames]{xcolor}
\usepackage[colorlinks=true, linkcolor=Maroon, urlcolor=Black]{hyperref}
\usepackage{amsthm}
\usepackage{amsfonts}
\usepackage{epsfig}
\usepackage{ textcomp }
\usepackage{mathtools}
\usepackage{authblk}
\usepackage{etoolbox}
\usepackage[bottom]{footmisc}
\usepackage{lipsum}
\usepackage{graphicx} 
\usepackage{enumitem}
 \newtheorem{thm}{Theorem}
 
 \newtheorem{prop}[thm]{Proposition}
 \newtheorem{lem}[thm]{Lemma}
 \newtheorem{prob}[thm]{Problem}

 \newcommand{\QED} {\hfill $\square$}

\DeclareMathOperator{\dom}{dom}
\DeclareMathOperator{\Area}{Area}
\DeclareMathOperator{\Vol}{Vol}
\DeclareMathOperator{\Span}{Span}
\title{Minimal and Maximal Distances in Metric Spaces}
\author{\v{Z}arko Ran\dj elovi\'c*}
\date{July 2025}

\begin{document}

\maketitle

\begin{abstract}
    Given functions $f,g: [n] \rightarrow [n]$, do there exist $n$ points $A_1,A_2,\ldots,A_n$ in some metric space such that $A_{f(i)},A_{g(i)}$ are the points closest and farthest from point $A_i$? In this paper we characterize precisely which pairs of functions have this property. Define $m(k)$ to be the maximal number such that any pair of functions $f,g:[m(k)]\rightarrow [m(k)]$ realizable in some metric space is also realizable in $\mathbb{R}^k$. We show that $m(k)$ grows exponentially in $k$. This answers a question of Croft. We also discuss what happens when looking at minimal and maximal distances separately.
    
\end{abstract}

\section{Introduction}
Let $n\ge 3$ be a positive integer and let $f,g:[n]\rightarrow [n]$ be two functions without fixed points, where as usual $[n]$ denotes the set $\{1,2,\ldots,n\}$. We will call a pair of such functions $(f,g)$ \textit{realizable in $(X,d)$} (where $(X,d)$ is a metric space) if there are $n$ points $A_1,\ldots,A_n\in X$ such that:
\begin{itemize}
    \item $d(A_i,A_j)$ are all distinct for different pairs $\{i,j\}$ with $i\neq j$.
    \item For each $i$, we have that $d(A_i,A_{f(i)})<d(A_i,A_j)$ for all $j\neq i,f(i)$.
    \item For each $i$, we have that $d(A_i,A_{g(i)})>d(A_i,A_j)$ for all $j\neq i,g(i)$.
\end{itemize}
Call a pair $(f,g)$ \textit{realizable} if it is realizable in any metric space $(X,d)$. We trivially observe that for any realizable pair $(f,g)$, we must have $f(i)\neq g(i)$ for each $i$. For convenience, we may view $f,g$ as directed graphs on the vertex set $[n]$. With this in mind, for any function $f:[n]\rightarrow [n]$, we define $G_f$ to be the directed graph on vertex set $[n]$ where there is an edge from $i$ to $j$ (which we will denote as $i\rightarrow j$) if and only if $f(i)=j$. We will now make a few definitions regarding directed graphs.

\makeatletter{\let\thefootnote\relax\footnote{*Mathematical Institute of the Serbian Academy of Sciences and Arts, Kneza Mihaila 36, Belgrade 11000, Serbia. Email: zarko.randjelovic@turing.mi.sanu.ac.rs}}

In a directed graph define the \emph{out-degree} of vertex $v$ to be the number of vertices $w$ such that there is an edge from $v\rightarrow w$. Similarly, define the \emph{in-degree} of a vertex $v$ to be the number of vertices $w$ such that $w\rightarrow v$ is an edge. We also define a \textit{source} in a directed graph to be a vertex of in-degree zero. The first result of this paper will describe precisely which pairs $(f,g)$ are realizable. 
\begin{thm}\label{genrealizable}
    Suppose that $n\ge 3$ is an integer and let $f,g:[n]\rightarrow [n]$ be two functions without fixed points such that for all $1\le i\le n$, we have $f(i)\neq g(i)$. Then the pair $(f,g)$ is realizable if and only if all of the following conditions hold:
    \begin{itemize}
        \item $G_f$ and $G_g$ do not have any cycles that are not of length two.
        \item $f\circ g$ has at most one fixed point.
        \item If $f(g(i))=i$, then $i$ is a source in $G_g$ and $g(i)$ is a source in $G_f$.
    \end{itemize}
\end{thm}

This result is surely folklore although we could not find it in the literature. The main aim of this paper is to see what happens if our metric space is some $\mathbb{R}^k$. Given a positive integer $n$, what is the minimal $k$ such that any realizable pair $(f,g)$ with $|\dom f|=n$ is also realizable in $\mathbb{R}^k$? This is a question of Croft \cite{Croft}. Notice that in our problem the only important aspect of the distances between the points is how they are ordered. Let $(f,g)$ be a realizable pair with $n=|\dom (f)|$ and suppose that $(f,g)$ is realized in $(X,d)$ with points $A_1,\ldots,A_n$. Then we can order the $2$-element subsets of $[n]$ in the following way:
$$\{i,j\}<\{s,l\}\ \  \text{if} \ \ d(A_i,A_j)<d(A_s,A_l).$$
We may pick real numbers $a_{\{i,j\}}$ to be all very close to $1$ (depending on $n$) and be in the same order as the $2$-sets. Now, it is not hard to show that we can pick points $B_1,\ldots,B_n$ in $\mathbb{R}^{n-1}$ such that $B_iB_j=a_{\{i,j\}}$ for any $i,j$, where $B_iB_j$ is the Euclidean distance between $B_i$ and $B_j$. This means that $(f,g)$ is in fact realizable in $\mathbb{R}^{n-1}$. The points $B_1,\ldots ,B_n$ are essentially an approximation of the vertices of a regular $(n-1)$-simplex.

Can we do it in less than $n-1$ dimensions? A natural question to ask is if we are given a positive integer $k$, what is the largest $n$ so that any realizable pair $(f,g)$ with $|\dom f|=n$ is in fact realizable in $\mathbb{R}^k$. We define $m(k)$ to be the largest such $n$. From our example above we see that $m(k)\ge k+1$. In Section \ref{section3}, we show that $m(k)$ actually grows exponentially. 
\begin{thm}\label{Rkrealizable}
    There are constants $A,B>0$ and $c,C>1$ such that for any positive integer $k$, we have that $Ac^k<m(k)<BC^k$.
\end{thm}
In the final section of this paper, we will show a few related results regarding minimal and maximal distances separately.
We now introduce some notation. Let $K_n$ be the complete undirected graph on the vertex set $[n]$. If $f$ is a function from $[n]$ to $[n]$ without fixed points, denote by $H_f$ the undirected graph on vertex set $[n]$ where $ij$ is an edge of $H_f$ if and only if either $f(i)=j$ or $f(j)=i$. We will always assume that $\mathbb{R}^k$ has the Euclidean metric and if $A,B\in \mathbb{R}^k$ for some $k$, we will let $AB$ denote the Euclidean distance between $A$ and $B$. Throughout this paper we will use standard notation for graphs and their parameters. See Bollob\'{a}s \cite{10.5555/7228} for general background.
\section{Proof of Theorem \ref{genrealizable}}
In order to prove Theorem \ref{genrealizable}, we will first prove a few lemmas that will give us that the three conditions are indeed necessary. Then, for any pair $(f,g)$ satisfying those conditions, we will give an example that will show that $(f,g)$ is realizable. For convenience, we say that a pair $(f,g)$ of functions from $[n]$ to $[n]$ is \emph{nice} if all of the following conditions hold:
\begin{itemize}
    \item For any $1\le i \le n$, we have that $i,f(i),g(i)$ are distinct.
    \item $G_f$ and $G_g$ do not have any cycles that are not of length two.
    \item $f\circ g$ has at most one fixed point.
    \item If $f(g(i))=i$, then $i$ is a source in $G_g$ and $g(i)$ is a source in $G_f$.
\end{itemize}
The following lemma shows that the first condition in Theorem \ref{genrealizable} is necessary.
\begin{lem}\label{no2cycles}
    Suppose that $(f,g)$ is a realizable pair. Then $G_f$ and $G_g$ do not have any cycles that are not of length two.
\end{lem}
\textit{Proof.} Let $n=|\dom f|$. Since $(f,g)$ is realizable, there is a metric space $(X,d)$ with $n$ points $A_1,\ldots A_n$ such that $d(A_i,A_j)$ are all distinct for different pairs $\{i,j\}$ with $i\neq j$, and for any distinct $1\le i,j\le n$, we have that $d(A_i,A_{f(i)})\le d(A_i,A_j)$ and $d(A_i,A_{g(i)})\ge d(A_i,A_j)$. Suppose that $G_f$ has a cycle $i_1\rightarrow i_2\rightarrow i_3\ldots i_k\rightarrow i_1$ where $k\ge 3$. Take $1\le l\le k$ to be such that $$d(A_{i_l},A_{i_{l+1}})=\min_{1\le j\le k} d(A_{i_j},A_{i_{j+1}})$$
where $i_{k+1}=i_1$. Without loss of generality assume that $l=1$. Then $d(A_{i_1},A_{i_2})<d(A_{i_2},A_{i_3})$. But since $i_2\rightarrow i_3$ is an edge in $G_f$, we have that $f(i_2)=i_3$. This gives a contradiction. We similarly prove that $G_g$ does not have any cycles that are not of length two. This completes the proof. \QED

In the next lemma we show that the remaining two conditions in Theorem \ref{genrealizable} are also necessary.
\begin{lem}\label{onefixedpoint}
    Suppose that $(f,g)$ is a realizable pair. Then we must have that $f\circ g$ has at most one fixed point. Moreover, if $i$ is that fixed point, then $i$ is a source in $G_g$ and $g(i)$ is a source in $G_f$.
\end{lem}
\textit{Proof.} Let $n=|\dom f|$. Since $(f,g)$ is realizable, there is a metric space $(X,d)$ with $n$ points $A_1,\ldots A_n$ such that $d(A_i,A_j)$ are all distinct for different pairs $i\neq j$, and for any distinct $1\le i,j\le n$, we have that \begin{equation}\label{7.1}
d(A_i,A_{f(i)})\le d(A_i,A_j)\ \ \text{and}\ \ d(A_i,A_{g(i)})\ge d(A_i,A_j).\end{equation}
First suppose that $i$ is a fixed point of $f\circ g$ and suppose that $i$ is not a source of $G_g$. That means that there is some $j$ such that $g(j)=i$. We first observe that $j\neq g(i)$. Indeed, if $j=g(i)$, then $f(j)=g(j)=i$ which is impossible. But now from (\ref{7.1}) we have that $$d(A_j,A_i)=d(A_j,A_{g(j)})>d(A_j,A_{g(i)})>d(A_{f(g(i))},A_{g(i)})=d(A_i,A_{g(i)})>d(A_i,A_j),$$ which gives a contradiction. This means that $i$ is a source in $G_g$. Similarly, we can prove that $g(i)$ is a source in $G_f$.\\
It remains to show that $f\circ g$ has at most one fixed point. We will also prove this by contradiction. Suppose that $i,j$ are two distinct fixed points of $f\circ g$. If $j=g(i)$, then $j$ is not a source of $G_g$ which is impossible. Similarly, $i\neq g(j)$. We also cannot have $g(i)=g(j)$ since if we did then $i=f(g(i))=f(g(j))=j$. Therefore, $i,g(i),j,g(j)$ are all distinct. But then using (\ref{7.1}) we have 
\begin{align*}
d(A_j,A_{g(i)})>d(A_{f(g(i))},A_{g(i)})=d(A_i,A_{g(i)})>d(A_i,A_{g(j)})>\\
>d(A_{f(g(j))},A_{g(j)})=d(A_j,A_{g(j)})>d(A_j,A_{g(i)})
\end{align*}
which is a contradiction. This completes the proof. \QED
\\
\\
We now prove the following lemma that shows us a bit about the structure of $G_f$.

 \begin{lem}\label{partitiontree}
     Suppose that $f:[n]\rightarrow [n]$ is a function without fixed points such that $G_f$ has no cycles of length more than two. Then there is some $m\ge 0$ and a partition $[n]=F_0\cup F_1\cup \ldots\cup F_m$ where $F_0$ contains precisely the points that are part of some $2$-cycle in $G_f$ and for every edge $i\rightarrow j$ of $G_f$, there is either some $1\le l\le m$ such that $i\in F_l,j\in F_{l-1}$ or $i,j\in F_0$.
 \end{lem}
 \textit{Proof.} Let $F_0=\{i\ |\ f(f(i))=i\}$. Since $f$ has no fixed points, we see that $F_0$ is precisely the set of points that belong to a $2$-cycle in $G_f$.  Define $F_1=\{i\ |\ i\not \in F_0, f(i)\in F_0\}$. We now iteratively define $F_k$ for integers $k\ge 2$. Let $$F_k=\{i\ |\ f(i)\in F_{k-1}\}.$$
We first show that the $F_k$ are disjoint. Trivially $F_0,F_1$ are disjoint. If $F_0,\ldots F_{k-1}$ are disjoint for some $k\ge 2$, then if $i\in F_k$. we must have that $f(i)\in F_{k-1}$. But then we cannot have that $i\in F_l$ for some $l<k$, because if $l>0$ then $f(i)\in F_{l-1}$ and if $l=0$ then $f(i)\in F_0$. Thus, $F_0,F_1,\ldots$ are disjoint sets by induction. Since $F_k\subset [n]$ for any $k$ there must be a largest $m$ such that $F_m\neq \emptyset$.\\
\\
We will now show that 
\begin{align}\label{7.2}
\cup_{j=0}^mF_j=[n].
\end{align}
Since for all $j$, we have that $F_j\subset [n]$, it follows that $\cup_{j=0}^mF_j\subset [n]$. \\
\\
 Let $i\in [n]$. We will show that $i\in \cup_{j=0}^mF_j$. For any $k\in \mathbb{Z}_{\ge 0}$, define $f^{(k)}(i)=f(f(\ldots f(i)\ldots ))$ where $f$ is applied $k$ times ($f^{(0)}(i)=i$). Now we consider the sequence $i,f(i),f^{(2)}(i),\ldots$. Since all of these numbers are in $[n]$, eventually some number in this sequence must repeat. Suppose that $k$ is the largest integer such that $i,f(i),f^{(2)}(i),\ldots f^{(k-1)}(i)$ are all distinct. Then there is a unique $j<k$ such that $f^{(j)}(i)=f^{(k)}(i)$. Since $f$ has no fixed points, we must have that $j<k-1$. But now $f^{(j)}(i)\rightarrow f^{(j+1)}(i)\rightarrow \ldots \rightarrow f^{(k-1)}(i)\rightarrow f^{(j)}(i)$ is a cycle in $G_f$ and therefore has length two. This means that $j=k-2$. If $k=2$, then we have that $f(f(i))=i$ so $i\in F_0$ and we are done. If $k>2$, then since $f(f(f^{(k-3)}(i)))\neq f^{(k-3)}(i)$, we have that $f^{(k-3)}(i)\not \in F_0$. But since $f(f^{(k-3)}(i))=f^{(k-2)}(i)\in F_0$ we have that $f^{(k-3)}(i)\in F_1$. \\
\\
Now we can show that $f^{(k-2-l)}(i)\in F_l$ by induction. We have just shown that it is true for $l=1$. Suppose that it is true for some  $1\le l\le k-3$. Since $f^{(k-2-l)}(i)\in F_l$, by definition of $F_{l+1}$ we have that $f^{(k-2-(l+1))}(i)\in F_{l+1}$. Thus, by induction, $f^{(k-2-l)}(i)\in F_l$ for any $1\le l\le k-2$. This means that $i\in F_{k-2}$ and hence $i\in \cup_{j=0}^mF_j$. This shows that $[n]\subset \cup_{j=0}^mF_j$ which then gives (\ref{7.2}).\\
\\
Notice that if $i\in F_0$ then $f(i)\in F_0$. But also, if $i\in F_l$ for some $l>0$, then $f(i)\in F_{l-1}$ so the partition satisfies the desired property which completes the proof. \QED
\\
\\
 For any pair $(f,g)$ of functions from the set $[n]$ to itself without fixed points, define \linebreak $T_{f,g}=H_f\cap H_g$ and define $T_f,T_g$ to be the graphs with $V(T_f)=V(T_g)=[n]$ and $E(T_f)=E(H_f)\backslash E(T_{f,g}),E(T_g)=E(H_g)\backslash E(T_{f,g})$.
 
\begin{lem}\label{Tfg}
    For any nice pair $(f,g)$, the graph $T_{f,g}$ has at most one edge and if it does have an edge, then that edge is precisely $ig(i)$, where $i$ is the unique fixed point of $f\circ g$.
\end{lem}
\textit{Proof.}  Since $G_f$ and $G_g$ do not have any common edges, the only way $H_f$ and $H_g$ could have a common edge is if there are $u,v\in [n]$ such that $f(u)=v$ and $g(v)=u$. We have that $f\circ g(v)=v$. Since $(f,g)$ is a nice pair, $v$ must be the unique fixed point of $f\circ g$ and the common edge of $H_f$ and $H_g$ is $vg(v)$.\QED 
\\
\\
\textit{Proof of Theorem \ref{genrealizable}.} We can see that Lemmas \ref{no2cycles},\ref{onefixedpoint} give that the three conditions in Theorem \ref{genrealizable} must hold for any realizable pair $(f,g)$. Now suppose that $f,g:[n]\rightarrow [n]$ are two functions without fixed points such that for all $i$, we have $f(i)\neq g(i)$ and the three conditions of Theorem \ref{genrealizable} are all satisfied. In other words, $(f,g)$ is a nice pair. We will show that the pair $(f,g)$ is realizable. \\
 
 We will  label the edges of $K_n$ in a way so that the order of the labels will correspond to the order of the distances in the required example. 
\\
\\
 Let $F_0,F_1,\ldots,F_m$ be the partition of $[n]$ given by Lemma \ref{partitiontree} for $G_f$ where $m\in \mathbb{Z}_{\ge 0}$. Similarly, let $G_0,\ldots,G_k$ be the partition of $[n]$ given by Lemma \ref{partitiontree} for $G_g$ where $k\in \mathbb{Z}_{\ge 0}$. For each edge $e$ of $T_f$, define the edge $e$ to be \emph{type} $0$ if both endpoints of $e$ are in $F_0$ and \emph{type} $i$ if one endpoint is in $F_i$ and the other is in $F_{i-1}$.
 We will order the edges $T_f$ in order $e_1,\ldots,e_{e(T_f)}$ from lowest to highest type (edges of the same type can be ordered arbitrarily). Likewise, we can define types of edges of $T_g$ and also order them in order $f_1,\ldots,f_{e(T_g)}$ in a similar fashion. 
\\
\\
We will now label the edges of $K_n$ as mentioned above. To do this we will construct a bijective map $c: E(K_n)\rightarrow [\binom{n}{2}]$ step by step. At any point we refer to edge $e$ as labeled if $c(e)$ has been defined up to that point and unlabeled otherwise. We first label the edges of $T_f$ and $T_g$. Define $c(e_i)=i$ for each $1\le i \le e(T_f)$ and define $c(f_i)=\binom{n}{2}-i+1$ for each $1\le i\le e(T_g)$. \\
\\
If $T_{f,g}$ has an edge, then by the proof of Lemma \ref{Tfg} there is a unique $q$ such that $f\circ g(q)=q$ and the only edge of $T_{f,g}$ is $qg(q)$. Notice that $qg(q)$ is unlabeled. Now suppose that there are $l$ so far unlabeled edges $x_1,\ldots,x_l$ incident to $q$ but not $g(q)$ where $l\in \mathbb{Z}_{\ge 0}$. Similarly, let $y_1,\ldots ,y_s$ be the so far unlabeled edges incident to $g(q)$ but not $q$, where $s\in \mathbb{Z}_{\ge 0}$. Define $c(x_j)=e(T_f)+j$ for all $1\le j\le l$ and define $c(y_j)=\binom{n}{2}+1-e(T_g)-j$ for all $1\le j\le s$. 
\\
\\
We now define $c$ arbitrarily on the remaining unlabeled edges of $K_n$. Consider now any $j\in [n]$. We will now show that 
\begin{align}\label{7.3}
    c(jf(j))<c(jr) \text{ for any } r\neq j,f(j).
\end{align}
First suppose that $j$ is not a fixed point of $g\circ f$. This means that $jf(j)\not \in E(H_g)$ and hence $jf(j)\in E(T_f)$. If $j\in F_0$, then by definition of $F_0$ in Lemma \ref{partitiontree} and by definition of $G_f$, we have that $jf(j)$ is the only edge incident to $j$ of type $0$. By definition of $c$ we have that (\ref{7.3}) holds in this case. If $j\in F_p$ for some $p\ge 1$, then by definition of $F_p$ from Lemma \ref{partitiontree} the only edge incident to $j$ of type at most $p$ is $jf(j)$. Thus, (\ref{7.3}) holds in this case as well. Finally, if $j$ is a fixed point of $g\circ f$, then $f(j)$ is the unique fixed point of $f\circ g$ and the unique edge of $T_{f,g}$ is $jf(j)$. We also have that $j$ is a source in $G_f$ and hence no edge in $T_f$ is incident to $j$. By definition of $c$ the edges incident to $j$ that are not edges of $T_f$ or $T_g$ have been labeled with values larger than $jf(j)$. Thus, (\ref{7.3}) holds in this case and hence (\ref{7.3}) always holds.
\\
\\
Similarly, we can show that for any $j\in [n]$ 
\begin{align}\label{7.4}
    c(jg(j))>c(jr) \text{ for any } r\neq j,g(j).
\end{align}
Now we are ready to construct our metric space. Let $1<a_1<a_2<\ldots a_{\binom{n}{2}}<2$ be arbitrary real numbers. Let $X=[n]$ and for any $i,j\in [n]$, define the distance 
\begin{align*}
d(i,j)=
  \bigg\{
  \begin{aligned}
  a_{c(ij)}, \text{ if } i\neq j \\
  0, \text{ if } i=j
  \end{aligned}.
\end{align*}
If $i,j,r\in [n]$, then $d(i,j)+d(j,r)\ge d(i,r)$ is clearly satisfied (if $i=j$ or $j=r$, we have equality and otherwise we have strict inequality). Thus, $(X,d)$ is indeed a metric space and the distances between all pairs of points are distinct. By (\ref{7.3}) and (\ref{7.4}), for each $i\neq j$, we have that $d(i,f(i))\le d(i,j)\le d(i,g(i))$. We also have that all distances $d(i,j)$ are distinct for different pairs $\{i,j\}$ with $i\neq j$. Thus, we have shown that the pair $(f,g)$ is realizable. \QED

\section{Proof of Theorem \ref{Rkrealizable}}\label{section3}
To prove Theorem \ref{Rkrealizable} we will first prove the upper bound for $m(k)$ which is easier. To do this we will need an estimate of the volume of a cap on a unit $(k-1)$-sphere. We then construct an example of a realizable pair, with the required number of points, that is not realizable in $\mathbb{R}^{k}$. We then move on to the lower bound where for any realizable pair $(f,g)$ with sufficiently low $|\dom f|$ (that is still exponential in $k$), we need to construct an example in $\mathbb{R}^k$ that realizes the pair $(f,g)$. To do this we will use the unit $(k-1)$-sphere in $\mathbb{R}^k$ and we will have that all the points in our example lie on that sphere. For any $k$, let $S^k$ be the unit $k$-sphere defined as $S^k=\{x\in \mathbb{R}^{k+1}\ |\ \lVert x\rVert=1\}$. All angles will be expressed in radians. The following will be a useful fact:
\begin{prop}\label{angle60}
     Suppose that $k,n\ge 3$ are positive integers and $A_1,\ldots,A_n\in \mathbb{R}^k$. If $A_iA_3<A_iA_j$ for any $i=1,2,j\neq i,3$, then $\angle A_1A_3A_2> \pi/3$.
\end{prop}
\textit{Proof.} We have that $A_1A_3<A_1A_2$ and $A_2A_3<A_2A_1$. Hence $A_1A_2$ is the longest side in $\triangle A_1A_2A_3$. Therefore, it is opposite of the largest angle and hence $\angle A_1A_3A_2> \pi/3$. \QED
\\
\\
Notice that if $v_1,v_2$ are unit vectors in $\mathbb{R}^k$, then the angle between them (that is less than $\pi$) is equal to $\angle (v_1,v_2)=\arccos{(v_1\cdot v_2)}$. Let $v_1,v_2,v_3$ be three unit vectors in $\mathbb{R}^k$. We can apply an orthogonal transformation, which is an angle-preserving map, to map them to three unit vectors in $\mathbb{R}^3$. We then have that $v_1,v_2,v_3\in S^2$. The spherical distance between two vectors in $S^2$ is equal to the angle between them. Hence, by the triangle inequality on the sphere, we have that 
\begin{align}\label{7.5}
    \angle (v_1,v_2)+\angle (v_2,v_3)\ge \angle (v_1,v_3).
\end{align}

For any vector $v\in S^{k-1}$ and angle $0<\phi<\pi$, define a \emph{hyperspherical cap with angle} $\phi$ to be the set of points $C_{v,\phi}=\{w\ |\ w\in S^{k-1}, \angle (w,v)\le \phi\}.$


To proceed we will need the surface area of a hyperspherical cap. Let $k\ge 2$ and suppose that we have the sphere $S^k$ in $\mathbb{R}^{k+1}$. Let $\phi\in [0,\pi/2]$. Li \cite{hyperspherical} showed that the area of the hyperspherical cap with angle $\phi$ is equal to 
\begin{align*}
    A_{k,\phi}=\frac{2\pi^{(k-1)/2}}{\Gamma\left(\frac{k-1}{2}\right)}\int_0^{\phi} \sin^{k-2} \theta \ d\theta
\end{align*}
where $\Gamma$ is the Gamma function. Thus, if $0<\phi_1<\phi_2\le \pi/2$ then 
\begin{align}\label{7.6}
    \frac{A_{k,\phi_1}}{A_{k,\phi_2}}=\frac{\int_0^{\phi_1} \sin^{k-2} \theta \ d\theta}{\int_0^{\phi_2} \sin^{k-2} \theta \ d\theta}.
\end{align}
Let $A_k$ denote the area of $S^k$. We have that $A_{k,\pi/2}=A_k/2$. We now prove the following lemma that will give us upper and lower exponential bounds on the ratio of areas of hyperspherical caps with fixed angles.

\begin{lem}\label{expbounds}
    Given $0<\phi_1<\phi_2\le \pi/2$ and any positive integer $k\ge 3$ we have that $$\frac{\phi_1\sin^{k-2}\left(\frac{\phi_1}{2}\right)}{2\phi_2\sin^{k-2}(\phi_2)}<\frac{A_{k,\phi_1}}{A_{k,\phi_2}}<\frac{2\phi_1\sin^{k-2}(\phi_1)}{(\phi_2-\phi_1)\sin^{k-2}\left(\frac{\phi_1+\phi_2}{2}\right)}.$$
\end{lem}

\textit{Proof. } To get the lower bound we can see that for any $0<\phi\le \pi/2$,
\begin{align}\label{7.7}
\int_0^{\phi} \sin^{k-2} \theta \ d\theta>\int_{\phi/2}^{\phi} \sin^{k-2} \theta \ d\theta>\frac{\phi}{2}\sin ^{k-2} (\phi/2)
\end{align}
and we also have that 
\begin{align}\label{7.8}
\int_0^{\phi} \sin^{k-2} \theta \ d\theta < \phi \sin^{k-2}(\phi).    
\end{align} 
Now, from (\ref{7.6}),(\ref{7.7}) for $\phi_1$ and (\ref{7.8}) for $\phi_2$ we have that $$\frac{\phi_1\sin^{k-2}\left(\frac{\phi_1}{2}\right)}{2\phi_2\sin^{k-2}(\phi_2)}<\frac{A_{k,\phi_1}}{A_{k,\phi_2}}.$$
To get the upper bound we can also see that 
$$\int_0^{\phi_2} \sin^{k-2} \theta \ d\theta>\int_{\frac{\phi_1+\phi_2}{2}}^{\phi_2} \sin^{k-2} \theta \ d\theta>\left( \frac{\phi_2-\phi_1}{2}\right)\sin^{k-2}\left(\frac{\phi_2+\phi_1}{2}\right).$$
Now using this and (\ref{7.8}) for $\phi_1$ we have that $$\frac{A_{k,\phi_1}}{A_{k,\phi_2}}<\frac{2\phi_1\sin^{k-2}(\phi_1)}{(\phi_2-\phi_1)\sin^{k-2}\left(\frac{\phi_1+\phi_2}{2}\right)}.$$
This completes the proof of the lemma. \QED
\\

We will now need some more lemmas that will be used to prove the lower bound in Theorem \ref{Rkrealizable}. To do this we will need to describe how to construct examples in $\mathbb{R}^k$ for realizable pairs $(f,g)$. In these constructions we will pick points one by one. The following lemma will be used to find a convenient order in which to construct the points.
\begin{lem}\label{pointsort}
     Let $T_1,T_2,T_3,T_4$ be forest graphs on the same vertex set $V$ and let $G=T_1\cup T_2\cup T_3\cup T_4$. Then there is a permutation $x_1,x_2\ldots x_{|V|}$ of the vertices of $G$ such that for any $1\le j\le |V|$, we have that $|\{l\ |\ l<j,x_jx_l\in E(G)\}|\le 7$.
\end{lem}

\textit{Proof.} Let $n=|V|$. We show that the lemma is true by induction on $n$. The base case $n=1$ is trivial. Suppose it is true for $n=s\ge 1$. Now suppose that $n=s+1$. Since $T_i$ are forests for $1\le i \le 4$, we have that $e(T_i)\le s$ for $1\le i\le 4$. Thus, the average degree of $G$ satisfies $$d=\frac{2e(G)}{s+1}\le \frac{8s}{s+1}<8.$$ Hence there is a vertex $y\in G$ with degree at most $7$. Let $T_i'$ be the induced subgraph of $T_i$ on the vertex set $V'=V\backslash \{y\}$ for $1\le i \le 4$ and let $G'=\cup_{i=1}^4 T_i'$. Now, by induction, we can pick an order $x_1,x_2,\ldots x_s$ of the vertices of $V'$ satisfying that for any $1\le j\le s$, we have that $|\{l\ |\ l<j,x_jx_l\in E(G')\}|\le 7$. Now, let $y=x_{s+1}$. We have that the order $x_1,x_2\ldots x_{s+1}$ of the vertices in $V$ satisfies the required conditions. \QED
\\

Our construction for the lower bound will be on the sphere $S^{k-1}$. For points $A,B\in S^{k-1}$, we have that $AB=2\sin \frac{\angle AOB}{2}$ where $O$ is the center of $S^{k-1}$. Thus, the larger $\angle AOB$ is the larger the distance. For a realizable pair $(f,g)$ with $\dom f=[n]$, define $H_{f,g}$ to be the graph on vertex set $[n]$ with edge set:\\

$\begin{aligned}
E(H_{f,g})=
  \bigg\{
  \begin{aligned}
  E(H_f\cup H_g)  , \text{ if } E(H_f\cap H_g)=\emptyset  \ \ \ \ \  \\
  E(H_f\cup H_g)\cup \{ik_1,jk_2|i,j\in [n],i\neq k_1,j\neq k_2\}, \text{ if } k_1k_2 \in E(H_f\cap H_g)
  \end{aligned}.
\end{aligned}$\\
\\
We now prove a lemma that partially describe how we will pick our angles in potential constructions in $\mathbb{R}^k$ for realizable pairs.
\begin{lem}\label{Hfg}
    Suppose that $(f,g)$ is a realizable pair with $\dom f=[n]$. If $E(H_f)\cap E(H_g)\neq \emptyset$, let $k_1,k_2$ be such that $f(k_1)=k_2$ and $g(k_2)=k_1$. Then we can pick $\alpha_{i,j}$ for all edges $ij\in E(H_{f,g})$ satisfying the following properties
    \begin{itemize}
        \item $\alpha_{i,j}=\alpha_{j,i}$ and all $\alpha_{i,j}$ for $i<j$ are distinct.
        \item For any edge $ij\in H_f\backslash H_g$, we have that $\arccos(1/500)<\alpha_{i,j}<\arccos(1/1000)$.
        \item For any edge $ij\in H_g\backslash H_f$, we have that $\arccos(-1/1000)<\alpha_{i,j}<\arccos(-1/500)$.
        \item For any edge $ij\in H_f\cup H_g$, we have that $\alpha _{i,f(i)}\le \alpha_{i,j}\le \alpha_{i,g(i)}$.
        \item If $E(H_f)\cap E(H_g)\neq \emptyset$, then for any $i\neq k_1,j\neq k_2$ such that $ik_1,jk_2\not \in E(H_f\cup H_g)$, we have that $\arccos (1/1000)< \alpha_{k_2,j}<\pi/2<\alpha_{k_1,i}<\arccos (-1/1000)$.
        \item If $E(H_f)\cap E(H_g)\neq \emptyset$ then $\alpha_{k_1,k_2}=\pi/2$.
    \end{itemize}
\end{lem}
\textit{Proof.} From the proof of Theorem \ref{genrealizable} it is clear that we can define an ordering with a bijection $c:E(K_n)\rightarrow \left[\binom{n}{2}\right]$ which satisfies that for any distinct $i,j\in [n]$, $c(if(i))\le c(ij)\le c(ig(i))$. Also, from that proof, we will have that the edges of $H_f\backslash H_g$ have the lowest values of $c$ (from $1$ to $|H_f\backslash H_g|$) and the edges of $H_g\backslash H_f$ have the highest values of $c$ (from $\binom{n}{2}-|H_g\backslash H_f|+1$ to $\binom{n}{2}$). If $E(H_f)\cap E(H_g)\neq \emptyset$, then we will also have $c(k_2j)<c(k_1k_2)<c(k_1i)$ for any $i\neq k_1,j\neq k_2$. Since $\arccos(1/500)<\arccos(1/1000)<\pi/2<\arccos(-1/1000)<\arccos(-1/500),$ we can clearly pick $\alpha_{i,j}$ to satisfy all the required properties.\QED
\\

For a hyperspherical cap $C_{v,\phi}\subset S^{k-1}$, define its \emph{boundary} to be $\partial C_{v,\phi}=\{w\ |\ w\in C_{v,\phi},\angle (w,v)=\phi\}$. Notice that $\partial C_{v,\phi}=\{w\ |\ w\in C_{v,\phi},w\cdot v=\cos \phi\}$. We will now show a useful lemma about vectors in $\mathbb{R}^k$.

\begin{lem}\label{7vecs}
    Let $0<\alpha<1/100$ be a real number. Suppose that $s\in [7]$ and $v_1,v_2,\ldots ,v_s\in \mathbb{R}^k$ are unit vectors, where $k\ge 7$ is an integer and $|v_i\cdot v_j|\le \alpha$ for any distinct $i,j$. Let $a_1,a_2,\ldots ,a_s$ be real numbers such that $|a_i|\le \alpha$ for all $i$. Then there is a unique vector $v\in \Span (v_1,\ldots ,v_s)$ such that $v\cdot v_i=a_i$ for all $i\in [s]$. Moreover, $\lVert v\rVert \le 8\alpha.$
\end{lem}
\textit{Proof.}  We first show that $\{v_i\ |\ 1\le i\le s\}$ is linearly independent. Suppose that there are $\mu_i$ not all zero for $1\le i\le s$ such that $\sum_{i=1}^s\mu_iv_i=0$. Without loss of generality let $|\mu_1|\ge |\mu_i|$ for all $2\le i\le s$. Then $$v_1=-\sum_{i=2}^s\frac{\mu_i}{\mu_1}v_i.$$
Taking the dot product with $v_1$ we have that 
$$1=|v_1\cdot v_1|=\left| -\sum_{i=2}^s\frac{\mu_i}{\mu_1}v_i\cdot v_1\right|\le \sum_{i=2}^s\left|\frac{\mu_i}{\mu_1}\right||v_1\cdot v_i|\le 6\alpha,$$
which is a contradiction. Thus, $v_1,\ldots ,v_s$ are linearly independent. By applying an appropriate orthogonal transformation we may assume without loss of generality that $\Span (v_1,\ldots ,v_s)=\mathbb{R}^s$. Now, if we consider the $s\times s$ matrix $A$ with row $i$ equal to the row vector $v_i$, then vectors $v\in \mathbb{R}^s$ satisfying $v\cdot v_i=a_i$ for all $i$ are precisely solutions to the equation 
$$Av=(a_1,a_2\ldots ,a_s)^{T}.$$
This has a unique solution since $v_1,\ldots ,v_s$ are linearly independent and therefore $A$ is invertible. Suppose that $v$ is that solution.

Now, let $v=\sum_{i=1}^s\lambda_iv_i$. Without loss of generality suppose that $|\lambda_1|\ge |\lambda_i|$ for all $2\le i\le s$. If $\lambda_1=0$ then $0=\lVert v\rVert<8\alpha$, so suppose that $\lambda_1\neq 0$. We have by the triangle inequality that 
\begin{align*}
    \alpha\ge |v\cdot v_1|=|\lambda_1|\left|v_1\cdot v_1 +\sum_{i=2}^s\frac{\lambda_i}{\lambda_1}v_i\cdot v_1\right|\ge \\ \ge |\lambda_1|\left(|v_1\cdot v_1|-\sum_{i=2}^s\left|\frac{\lambda_i}{\lambda_1}\right||v_i\cdot v_1|\right)\ge |\lambda_1|(1-6\alpha)
\end{align*}
since $v_1$ is a unit vector and $|v_1\cdot v_i|\le \alpha$ for $2\le i\le s$. Thus, $|\lambda_1|\le \frac{\alpha}{1-6\alpha}$ and hence $|\lambda_i|\le \frac{\alpha}{1-2\alpha}$ for all $i$. Now, by the triangle inequality and since $v_i$ are unit vectors, we have that 
$$\lVert v \rVert \le \sum_{i=1}^s|\lambda_i|\le \frac{7\alpha}{1-6\alpha}<8\alpha$$
since $\alpha<1/100$. This completes the proof of the lemma. \QED
\\

We will now show a lemma that tells us about the intersection of the boundaries of at most seven hyperspherical caps.
\begin{lem}\label{capintersection}
    Let $\alpha<1/200$ be a positive real number. Suppose that $2\le s\le 8$ and let $v_1,v_2,\ldots,v_{s}$ be unit vectors in $\mathbb{R}^k$ where $k\ge 9$ is an integer. Let $\alpha_i\in (0,\pi)$ for $i\in [s]$ be real numbers such that $|\cos (\alpha_i)|\le \alpha$ for all $i$ and $\cos \alpha_s\ge \alpha/2$. Suppose that $|v_i\cdot v_j|\le \alpha$ for all distinct $i,j\in [s]$. Let $C_{v_i,\alpha_i}$ for $i\in [s]$ be hyperspherical caps in $S^{k-1}$. We then have that $T=\cap_{i=1}^{s-1}\partial C_{v_i,\alpha_i}$ is a $(k-s)$-dimensional sphere of radius not smaller than $\sqrt{1-64\alpha^2}$. Moreover,
    $$\frac{\Area (C_{v_s,\alpha_s}\cap T)}{\Area (T)}\le \frac{A_{k-s,\beta}}{A_{k-s}}$$, where $\beta=\arccos (\frac{\alpha-128\alpha^2}{2})$
    and the areas are both $(k-s)$-dimensional areas.
\end{lem}
\textit{Proof.} Let $S=\Span (v_1,v_2,\ldots ,v_{s-1})$ and let $S^{\perp}$ be the orthogonal complement of $S$ in $\mathbb{R}^k$. Let $a_i=\cos \alpha_i$ for all $1\le i\le s$. By Lemma \ref{7vecs} there is a unique vector $v\in S$ such that $|v\cdot v_i|=a_i$ for $i\in [s-1]$ and we also have that $\lVert v \rVert \le 8\alpha<1$. From the proof of that lemma we also know that $v_1,v_2,\ldots v_{s-1}$ are linearly independent and hence $\dim S=s-1$. For any vector $w\in S^{\perp}$ with $\lVert w\rVert=\sqrt{1-\lVert v\rVert^2}$, we have that $\lVert v+w\rVert^2=\lVert v\rVert^2+\lVert w\rVert^2=1$. Thus, $v+w\in S^{k-1}$. We also have that $(v+w)\cdot v_i=a_i$ for all $i\in [s-1]$. On the other hand, suppose that $x\in S^{k-1}$ is such that $x\cdot v_i=a_i$ for all $i\in [s-1]$. Let $x=x_1+x_2$ be the decomposition of $x$ such that $x_1\in S,x_2\in S^{\perp}$. Then we have that $x_1\cdot v_i=a_i$ for $i\in [s-1]$ and hence by Lemma \ref{7vecs} $x_1=v$.
But then since $\lVert x_1\rVert^2+\lVert x_2\rVert^2=1$, we have that 
$\lVert x_2\rVert=\sqrt{1-\lVert v\rVert^2}$. Since $\partial C_{v_i,\alpha_i}=\{x\ |\ x\in S^{k-1},x\cdot v_i=a_i\}$, from above we have that 
$$T=\{x\ |\ x\in S^{k-1},x\cdot v_i=a_i \text { for all } i\in [s-1]\}=\{v+w\ |\ w\in S^{\perp},\lVert w\rVert=\sqrt{1-\lVert v\rVert^2}\}$$
is a $(k-s)$-dimensional sphere (since $\dim S^{\perp}=k-s+1$) of radius $\sqrt{1-\lVert v\rVert^2}\ge \sqrt{1-64\alpha^2}>0$. Thus, $T$ is non-empty. We now decompose the vector $v_s$ into $v_s=v_s'+w_s$, where $v_s'\in S$ and $w_s\in S^{\perp}$. We have that $v_s'\cdot v_i=v_s\cdot v_i$ for $i\in [7]$. Notice that by Lemma \ref{7vecs} we have that $v_s'$ is the unique vector in $S$ satisfying $v_s'\cdot v_i=v_s\cdot v_i$ for $i\in [s-1]$ and hence $\lVert v_s'\rVert\le 8\alpha$. 

Let $M=\{w\ |\ w\in S^{\perp},\lVert w\rVert^2=1-\lVert v\rVert^2\}$. Applying an appropriate orthogonal transformation we may assume that $S^{\perp}=\mathbb{R}^{k-s+1}$. Now we have that $M=\sqrt{1-\lVert v\rVert^2}S^{k-s}$. For any $w\in M$, we see that $v+w\in C_{v_s,\alpha_s}$ if and only if $(v+w)\cdot v_s\ge a_s$. But we have that 
$$(v+w)\cdot v_s=(v+w)\cdot (v_s'+w_s)=v\cdot v_s'+w\cdot w_s.$$
Thus, $v+w\in C_{v_s,\alpha_s}$ if and only if $w\cdot w_s\ge a_s-v\cdot v_s'$. Since $\lVert v_s'\rVert\le 8\alpha$, we know that $w,w_s$ are non zero and hence $v+w\in C_{v_s,\alpha_s}$ if and only if

\begin{align}\label{7.9}
\frac{w}{\lVert w\rVert}\cdot \frac{w_s}{\lVert w_s\rVert}\ge \frac{a_s-v\cdot v_s'}{\lVert w\rVert \lVert w_s\rVert}.
\end{align}
We know that $\lVert w_s\rVert,\lVert w\rVert\le 1$ and $|v\cdot v_s'|\le \lVert v\rVert \lVert v_s'\rVert\le 64\alpha^2$. Therefore, $a_s-v\cdot v_s\ge \alpha/2-64\alpha^2>0$ and hence 
$$\frac{a_s-v\cdot v_s'}{\lVert w\rVert \lVert w_s\rVert}\ge \frac{\alpha-128\alpha^2}{2}=\cos \beta>0.$$

Since $\frac{w}{\lVert w\rVert},\frac{w_s}{\lVert w_s\rVert}\in S^{k-s}$, we have that if (\ref{7.9}) is satisfied, then $\frac{w}{\lVert w\rVert}\in C_{\frac{w_s}{\lVert w_s\rVert},\beta}^{k-s}$, where the superscript indicates that the cap is in $S^{k-s}$. Since $M=T-v$ and $\beta<\pi/2$, we have that 
$$\frac{\Area (C_{v_s,\alpha_s}\cap T)}{\Area (T)}\le \frac{A_{k-s,\beta}}{A_{k-s}}$$ 
which completes the proof of the lemma. \QED \\

We are now ready to prove Theorem \ref{Rkrealizable}. First we prove the easier upper bound.\\
\\
\textit{Proof of upper bound.} Suppose that $k\ge 4$. Let $\epsilon=\sin (\pi/12)$ and let $\alpha=\frac{1}{6\sin^3(\pi/12)}$. By Lemma \ref{expbounds} we have that 
\begin{align}\label{7.10}
    \frac{A_{k-1,\pi/6}}{A_{k-1}}=\frac{A_{k-1,\pi/6}}{2A_{k-1,\pi/2}}>\frac{\sin^{k-3}(\pi/12)}{6}=\alpha \epsilon^k.
\end{align}
Now suppose that $m(k)>1+\frac{1}{\alpha \epsilon^k}$. Let $n=m(k)$. We know that $m(k)\ge 3$. Consider the following pair of functions $f,g:[n]\rightarrow [n]$ given by 

\begin{align*}
f(i)=
  \bigg\{
  \begin{aligned}
  n, \text{ if } i<n \\
  1, \text{ if } i=n
  \end{aligned}
\end{align*}
\begin{align*}
g(i)=
  \bigg\{
  \begin{aligned}
  1, \text{ if } 1<i<n \\
  2, \text{ if } i=1,n \ \
  \end{aligned}.
\end{align*}

We can see that by Theorem \ref{genrealizable} $(f,g)$ is a realizable pair. Since $n=m(k)$, there are distinct points $X_1,\ldots X_n\in \mathbb{R}^k$ such that for every $i$, we have that $X_iX_{f(i)}<X_iX_j$ for any $j\neq i,f(i)$. Let $v_i=\overrightarrow{X_nX_i}$ for each $1\le i\le n-1$. By Proposition \ref{angle60} we have that $\angle X_iX_nX_j> \pi/3$ for any $1\le i<j<n$ and hence $\angle (v_i,v_j)>\pi/3$. Let $w_i=\frac{v_i}{||v_i||}$. Now, by (\ref{7.5}) this means that the hyperspherical caps $C_{w_1,\pi/6},C_{w_2,\pi/6}\ldots C_{w_{n-1},\pi/6}$ in $S^{k-1}$ are pairwise disjoint. Therefore, we have that $(n-1)A_{k-1,\pi/6}\le A_{k-1}$. Thus, by (\ref{7.10}) we have that $$A_{k-1}\ge (n-1)A_{k-1,\pi/6}>\frac{A_{k-1,\pi/6}}{\alpha \epsilon^k}>A_{k-1},$$
which is a contradiction. Therefore, we have that $m(k)\le 1+\frac{1}{\alpha \epsilon^k}$. Let $C=\frac{2}{\epsilon}$. We have that $\epsilon<1$ so $C>1$ and there is a large enough $D$ such that for all $k\ge D$, we have that $2^{k}\alpha>1$ and $\frac{1}{\alpha \epsilon^k}>1$. Thus, for all $k\ge D$, we have that 
$$m(k)\le \frac{2}{\alpha \epsilon^k}\le \frac{2^{k+1}}{\epsilon^k}=2C^k.$$
Now, there is certainly a large enough $B$ so that for all positive integers $k$, we have that $m(k)\le BC^k$, which completes the proof. \QED

We will now prove the lower bound.

\textit{Proof of lower bound.} Let $\alpha=1/500$ and let $\beta=\arccos\left(\frac{\alpha-128\alpha^2}{2}\right)$. It is easy to see that $m(k)\ge 3$ for all $k$. We may assume that $k\ge 10$, since we can always pick $A>0$ to be small enough so that the lower bound holds for $1\le k\le 9$ as well. Suppose that $(f,g)$ is a realizable pair of functions with $\dom f=[n]$ and $n<Ac^k+1$, where $$c=\frac{\sin \left( \frac{\pi/2+\beta}{2}\right)}{\sin \beta} \text { and } A=\frac{(\pi/2-\beta)\sin^{10} \beta}{2\beta \sin ^{10}\left(\frac{\pi/2+\beta}{2}\right)}.$$ Since $\alpha=1/500$, we have that $\beta<\pi/2$ and hence $c>1$ and $A>0$. Let $\alpha_{i,j}$ for all $ij\in E(H_{f,g})$ satisfy the properties in Lemma \ref{Hfg}. Let $\delta=\arccos (\alpha/2)$.

If $H_f$ had a cycle say $i_1i_2\ldots i_k$ with $k\ge 3$, we may assume that $i_1\rightarrow i_2$ is an edge of $G_f$. But then since all vertices of $G_f$ have out-degree one and $i_1i_k\in E(H_f)$, we have that $i_k\rightarrow i_1$ is an edge of $G_f$. If $i_j\rightarrow i_{j+1}$ is an edge of $G_f$ for $j\ge 2$ ($i_{k+1}=i_1$), then we have that $i_{j-1}\rightarrow i_j$ is an edge of $G_f$. Thus, by induction, $i_1\rightarrow i_2\rightarrow \ldots \rightarrow i_k$ is a cycle in $G_f$, which is impossible. Therefore, $H_f$ has no cycles. Similarly, $H_g$ has no cycles.

 Let $T_1=H_f,T_2=H_g$. If $E(H_f\cap H_g)=\emptyset$, let $T_3=T_4=E_n$, where $E_n$ is the graph on vertex set $[n]$ with no edges. Otherwise, by Lemma \ref{Tfg} there is exactly one edge $k_1k_2$ in $H_f\cap H_g$ where $f(k_1)=k_2$ and $g(k_2)=k_1$. Suppose that $W_i$ is the graph on vertex set $[n]$ with edge set $\{jk_i\ |\ j\neq k_i\}$ for $i=1,2$. Notice that $$H_{f,g}=T_1\cup T_2\cup T_3\cup T_4.$$ We now apply Lemma \ref{pointsort} on $T_1,T_2,T_3,T_4$ to get a permutation $x_1,x_2\ldots x_n$ of $[n]$ such that for any $1\le j\le n$, we have that $|\{l\ |\ l<j, x_jx_l\in E(H_{f,g})\}|\le 7$. We may assume without loss of generality that $x_i=i$. Let $O$ be the coordinate origin of $\mathbb{R}^k$. We are now going to construct points $X_1,X_2,\ldots X_n\in S^{k-1}$ that satisfy the following properties:
\begin{itemize}
    \item $\delta< \angle X_iOX_j< \pi-\delta$ for any distinct $i,j\in [n]$ such that $ij\not \in E(H_{f,g}),$
    \item $\angle X_iOX_j=\alpha_{i,j}$ if $ij\in E(H_{f,g}).$ 
\end{itemize}
To do this we will construct these points inductively. First we pick $X_1$ to be an arbitrary point on $S^{k-1}$. Now suppose that we have so far constructed $X_1,X_2,\ldots X_{s-1}$ for some $2\le s\le n$ and that they satisfy the required properties. Let $v_i=\overrightarrow{OX_i}$ for $1\le i\le s-1$. We will now construct $X_s$. We know that there are at most seven edges $is$ in $H_{f,g}$ with $i<s$. In the remainder of the proof all hyperspherical caps are taken to be in $S^{k-1}$. 

 In order for the second property to hold it is enough that for any $is\in E(H_{f,g})$ with $1\le i\le s-1$, we have $X_s\in \partial C_{v_i,\alpha_{i,s}}$. Thus, we need $X_s$ to be in the intersection of some $t$ cap boundaries where $0\le t\le 7$. Assume that those cap boundaries are $\partial C_{w_i,\beta_i}$ for $1\le i\le t$, where the vectors $w_i$ for $1\le i\le t$ are some of the vectors $v_1,v_2,\ldots v_{s-1}$ and $\beta_i$ for $1\le i\le t$ are some of the $\alpha_{j,s}$. Since $\alpha =1/500$, by Lemma $\ref{Hfg}$ we know that for any distinct $1\le i,j\le t$, we have that $|w_i\cdot w_j|\le \alpha$.

On the other hand, if $is\not \in E(H_{f,g})$, then for the first property, it is enough that $X_s\not \in C_{v_i,\delta}$ and $X_s\not \in C_{-v_i,\delta}$. We know that $X_s$ needs to avoid at most $2(s-1)$ caps with angle $\delta$. Let those caps be 
$$C_{r_1,\delta},c_{r_2,\delta},\ldots C_{r_l,\delta},$$
where $0\le l\le 2(s-1)$ and every vector $r_i$ is one of the vectors $\pm v_1,\pm v_2,\ldots ,\pm v_{s-1}$. If $t\ge 1$ let 
$$T=\cap_{i=1}^t\partial C_{w_i,\beta_i}$$
and if $t=0$ let $T=S^{k-1}$. We will show that we can select $X_s\in T$ such that for every $1\le i\le l$, we also have that $X_s\not \in C_{r_i,\delta}$. 

First suppose that $t\ge 1$. We know that $r_i=\pm v_{i'}$ for some $i'$ such that $ii'\not \in E(H_{f,g})$. Hence none of the $w_1,\ldots w_t$ are equal to $v_{i'}$ but they are all among $v_1,v_2\ldots v_{s-1}$. Therefore, $|r_i\cdot w_j|\le \alpha$ for all $j\in [t]$. Thus, by Lemma \ref{capintersection} we have that for any $1\le i\le l$,
$$\frac{\Area (T\cap C_{r_i,\delta})}{\Area (T)}\le \frac{A_{k-t-1,\beta}}{A_{k-t-1}}.$$ 
If $t=0$, then we still have
$$\frac{\Area (T\cap C_{r_i,\delta})}{\Area (T)}\le \frac{A_{k-1,\delta}}{A_{k-1}}\le \frac{A_{k-t-1,\beta}}{A_{k-t-1}}$$ 
since $\beta>\delta$. Either way by Lemma \ref{expbounds}
$$\frac{A_{k-t-1,\beta}}{A_{k-t-1}}=\frac{A_{k-t-1,\beta}}{2A_{k-t-1,\pi/2}}\le \frac{\beta \sin^{k-t-3}\beta }{(\pi/2-\beta)\sin^{k-t-3} (\frac{\pi/2+\beta}{2})}\le \frac{1}{2Ac^k}$$
because $t\le 7$ and $\beta<\pi/2$. Thus, by the union bound
$$\frac{\Area (T\cap (\cup_{i=1}^l C_{r_i,\delta}))}{\Area (T)}\le \sum_{i=1}^l\frac{\Area (T\cap C_{r_i,\delta})}{\Area (T)}\le l\frac{1}{2Ac^k}\le \frac{2(s-1)}{2Ac^k}\le \frac{2(n-1)}{2Ac^k}<1$$ 
and hence there is a point $X_k\in T\backslash (\cup_{i=1}^l C_{r_i,\delta})$.
Thus, we ensure that the two required properties for angles $\angle X_iOX_j$ hold and by induction we can pick points $X_1,X_2\ldots X_n$ so that the two required properties hold. 

Now suppose that $i,j$ are distinct numbers in $[n]$ with $j\neq f(i)$. If $if(i)$ is not an edge in $H_g$, then it is an edge in $H_f\backslash H_g$. Hence, by definition of $\alpha_{i,j}$ from Lemma \ref{Hfg}, we have that if $ij\in E(H_f\cup H_g)$, by the fourth and first property of Lemma \ref{Hfg}$\angle X_iOX_{f(i)}<\angle X_iOX_j$. But also, if $if\not \in E(H_f\cup H_g)$, by the second property we have that $\angle X_iOX_{f(i)}<\delta\le \angle X_iOX_j$. On the other hand, if $if(i)\in H_g$ then $E(H_f\cap H_g)\neq \emptyset$, so let $k_1,k_2$ be defined as in Lemma \ref{Hfg}. We have that $i=k_1,f(i)=k_2$ and hence from the fourth, fifth and sixth property in Lemma \ref{Hfg} we have that $\angle X_iOX_{f(i)}=\pi/2<\alpha_{k_1,j}=\angle X_iOX_j$. Thus, for any distinct $i,j\in [n]$ with $j\neq f(i)$, we have that $X_iX_{f(i)}=2\sin \frac{\angle X_iOX_{f(i)}}{2}<2\sin  \frac{\angle X_iOX_j}{2}=X_iX_j$. Similarly, for any distinct $i,j\in [n]$ with $j\neq g(i)$, we have that $X_iX_{g(i)}<X_iX_j$. 

Now we just need to ensure that $X_iX_j$ are all distinct. Since the inequalities $X_iX_{f(i)}<X_iX_j$ for all $i,j\neq i,f(i)$ and $X_iX_{g(i)}>X_iX_j$ for all $i,j\neq i,g(i)$ are all strict, for any sufficiently small perturbation of points $X_1,X_2,\ldots X_n$, those inequalities will still hold. It is easy to see that we can find a small perturbation such that $X_iX_j$ are all distinct and all the inequalities still hold (we can move points one by one and we just need to avoid finitely many spheres and hyperplanes which have $k$-dimensional volume $0$). The perturbed points will not necessarily lie on $S^{k-1}$. This shows that $(f,g)$ is realizable in $\mathbb{R}^k$. Since there is an integer in $[Ac^k,Ac^k+1)$, this completes the proof of Theorem \ref{Rkrealizable}. \QED
\section{Related Results}

In this section we will consider minimal and maximal distance functions separately. It turns out that maximal distances are not so interesting. If $(X,d)$ is a metric space, we will define a function $g:[n]\rightarrow [n]$ without fixed points to be \emph{max-realizable in }$(X,d)$ if there are points $A_1,A_2\ldots A_n\in X$ such that all distances $d(A_i,A_j)$ are distinct for all different pairs $\{i,j\}$ with $i\neq j$ and for any distinct $i,j\in [n]$, we have that $d(A_i,A_{g(i)})\ge d(A_i,A_j)$. Define $g$ to be \emph{max-realizable} if it is max-realizable in some metric space.

\begin{thm}\label{maxrealizable}
    Any function that is max-realizable is also max-realizable in $\mathbb{R}^2.$
\end{thm}
Our construction will involve ellipses and we will need a few lemmas. The first lemma will be useful when choosing points on the ellipses.
\begin{lem}\label{ellipsedistance}
Let $C$ be the ellipse given by the equation $x^2+4y^2=4$. Suppose that $b$ is any real number in $(0,1/3)$. Let $P_b=(-2\sqrt{1-b^2},-b)$ and for any $y\in (0,1)$, let $Q_y=(2\sqrt{1-y^2},y)$. Define the function $g_b:(0,1)\rightarrow \mathbb{R}$ to be $g_b(y)=P_bQ_y^2$. Then there is an $m_b\in (0,1)$ (depending on $b$) such that $g_b(y)$ is strictly increasing on $(0,m_b]$ and strictly decreasing on $[m_b,1)$. Moreover, if $f:(0,1/3)\rightarrow (0,1)$ is defined as $f(b)=m_b$, then $f$ is a strictly increasing continuous function and $Q_{m_b}$ is the unique farthest point on $C$ from $P_b$.
\end{lem}
\textit{Proof.} Notice that $P_b,Q_y\in C$. We have that $$g_b(y)=4(\sqrt{1-y^2}+\sqrt{1-b^2})^2+(y+b)^2.$$ Computing the derivative with respect to $y$ we have that 
\begin{align}\label{7.11}
    g_b'(y)=-6y-8y\frac{\sqrt{1-b^2}}{\sqrt{1-y^2}}+2b.
\end{align}
We will show that there is a unique $m_b\in (0,1)$ such that $g'(m_b)=0$. To have $g_b'(y)=0$, dividing by $2$ and rearranging (\ref{7.11}), we need 
\begin{align}\label{7.12}
    3y-b=\frac{-4y\sqrt{1-b^2}}{\sqrt{1-y^2}}.
\end{align}
Since the right-hand side is negative, we must have $y<b/3$. We can also see that for $y\ge b/3$, we have that $g_b'(y)<0$. Now consider (\ref{7.12}) with $y\in (0,b/3)$. We have that 
\begin{align}\label{7.13}
    \frac{(3y-b)^2(1-y^2)}{y^2}=16(1-b^2).
\end{align}
Define $h_b(y)$ to be the left-hand side of (\ref{7.13}). We have that \begin{align*}
h_b'(y)=-18y+6b+6\frac{b}{y^2}-2\frac{b^2}{y^3}=\left(6-\frac{2b}{y^3}\right)(b-3y).
\end{align*}
We know that $1>b>3y$ and hence also $\frac{2b}{y^3}>\frac{2\cdot 27b}{b^3}>6$. This means that $h_b'(y)<0$ and hence $h_b(y)$ is decreasing on $(0,b/3)$. Since $h_b(b/3)=0$ and as $y\rightarrow 0$ $h_b(y)\rightarrow \infty$, we have that $(7.13)$ has a unique solution on $(0,b/3)$. This means that $g_b'(y)$ has a unique zero $m_b$ in $(0,1)$ which lies in $(0,b/3)$. Since $g_b'(y)$ is continuous and as $y\rightarrow 0$ $g_b'(y)\rightarrow 2b>0$, we have that $g_b'(y)$ is positive on $(0,m_b)$ and negative on $(m_b,1)$. So $g_b(y)$ is strictly increasing on $(0,m_b]$, strictly decreasing on $[m_b,1)$ and for all $y\in (0,1)\backslash \{m_b\}$, we have that $g_b(m_b)>g_b(y)$.\\

Now we will show that $f$ is continuous. Suppose that $b_n$ is a sequence of reals in $(0,1/3)$ such that $b_n\rightarrow b$ as $n\rightarrow \infty$ for some $b\in (0,1/3)$. We will show that $m_{b_n}\rightarrow m_b$. Assume the contrary. Then there is an $\epsilon>0$ such that there is an infinite subsequence $b_{i_n}$ such that $|m_{b_{i_n}}-m_b|>\epsilon$ for all $n$. By moving to a convergent subsequence, by compactness, we may assume that $m_{b_{i_n}}\rightarrow a$ for some $a$ such that $|a-m_b|\ge \epsilon$ and $a\ge 0$. Now, from (\ref{7.13}) we see that $h_{b_{i_n}}(m_{b_{i_n}})\rightarrow 16(1-b^2)=h_b(m_b)$, but also by continuity $h_{b_{i_n}}(m_{b_{i_n}})\rightarrow h_b(a)$ if $a>0$, while $h_{b_{i_n}}(m_{b_{i_n}})\rightarrow \infty$ if $a=0$. Thus, $a>0$ and since $m_{b_{i_n}}<b_{i_n}/3$, we must have that $a=\lim m_{b_{i_n}}\le \lim b_{i_n}/3= b/3$. Since $h_b(b/3)=0$, we have that $a<b/3$, but then we have that $a=m_b$, which is a contradiction. Thus, we have shown that $f$ is indeed continuous.\\

Now suppose that $b<c$. From (\ref{7.13}) and since $m_b<b/3<c/3$, we know that $h_c(m_b)>h_b(m_b)=16(1-b^2)>16(1-c^2)$. Since $h_c$ is decreasing on $(0,c/3)$, we have that $m_c>m_b$ and therefore $f$ is strictly increasing. Now we have from above that for any $b<1/3$, the farthest point to $P_b$ on $C$ (where $C$ refers only to the curve and not the area inside) in the open upper right quadrant is indeed $Q_{m_b}$. Since the tangents to the ellipse at points $A=(2,0)$ and $B=(0,1)$ are perpendicular to the $x,y$ axes respectively, we know that there are points $E,F$ on $C$ in the open upper right quadrant such that $\angle P_bAE,\angle P_bBF$ are obtuse. This means that neither $A$ nor $B$ can be the farthest point from $P_b$ on $C$. Now, for any point $(p,q)\in C$, if we consider the points $(p,q),(p,-q),(-p,q),(-p,-q)$, then the point (that appear twice if either coordinate is zero)  that is in the closed upper right quadrant has the highest distance to $P_b$ in both the $x,y$ coordinates. This means that $Q_{m_b}$ is indeed the unique farthest point from $P_b$ on $C$. \QED

In order to show that any function $g:[n]\rightarrow [n]$ is max-realizable in $\mathbb{R}^2$ we will first split $G_g$ into connected components and show that they are all max-realizable in certain types of subsets of $\mathbb{R}^2$. Then we will combine the examples for each component to show that $f$ is max-realizable in $\mathbb{R}^2$. For any $\delta>0$ and $X\in \mathbb{R}^k$, let $B_{\delta}(X)=\{P\in \mathbb{R}^k\ |\ XP<\delta \}$. The following lemma will deal with one component.
\begin{lem}\label{onecomponentrealizable}
    Suppose that $g:[n]\rightarrow [n]$ is a function without fixed points such that $H_g$ is connected and $G_g$ has no cycles of length bigger than two. Then $g$ is max-realizable in $\mathbb{R}^2$. Moreover, for any distinct points $A,B\in \mathbb{R}^2$ and any $\epsilon>0$, there is an example of points $X_1,X_2,\ldots X_n$ that realizes $g$ such that all the points are in $B_{\epsilon}(A)\cup B_{\epsilon}(B)$.
\end{lem}
\textit{Proof.} By rotation, translation and scaling we may assume that $A=(-2,0),B=(2,0)$. As above let $C$ be the ellipse given by the equation $x^2+4y^2=4$. We know that $A,B\in C$. We will show that we can pick the points $X_1,X_2,\ldots X_n$ to be on $C$. By Lemma \ref{partitiontree} there is a non-negative integer $m$ and a partition $[n]=F_0\cup F_1\cup \ldots \cup F_m$ with the property that if $i\rightarrow j$ is an edge of $G_g$, we have that either $i,j\in F_0$ or there is some $s>0$ such that $i\in F_s,j\in F_{s-1}$. From the proof of Lemma \ref{partitiontree} we have that $G_g[F_0]$ is precisely the set of points belonging to a $2$-cycle in $G_g$. We also know that since the out-degree of any vertex in $G_g$ is one, we have that $G_g[F_0]$ is a union of disjoint $2$-cycles and from any vertex $i\in [n]$ there is a directed path to exactly one of these cycles. Suppose that $C_1,C_2,\ldots C_s$ are all of the $2$-cycles in $F_0$. Then for $1\le i\le s$, let 
$$H_i=\{j\in [n]\ |\text{ there is a path from }j\text{ to }C_i\}.$$
By above $H_i$ are disjoint. If there is an edge say $u\rightarrow v$ from $H_i$ to $H_j$ for any $i\neq j$, then we would have $u\in H_j$, which is a contraction. This means that $H_g[H_1],H_g[H_2],\ldots H_g[H_s]$ are all disconnected from each other. Thus, $s=1$ and $G_g[F_0]=C_1$. \\

Let $C_1=\{a_0,b_0\}$. We will define the set $$S_A=\{i\in [n]\ |\text{ the shortest path in }G_g\text{ from }i\text{ to }C_1\text{ ends at }a_0\}.$$ Similarly, define $S_B$. For $1\le i\le m$, let $A_i=F_i\cap S_A$ and let $B_i=F_i\cap S_B$. Let $f$ be defined as in Lemma \ref{ellipsedistance}.
Let $X_{a_0}=A,X_{b_0}=B$. For any $b,y\in (0,1/3)$, define $P_b,Q_y$ as in Lemma \ref{ellipsedistance}. Let $0<b_m<1/3$ be such that $P_{b_m}A<\epsilon$ and define $b_{m-1},b_{m-2},\ldots b_1$ such that $f(b_i)=b_{i-1}$ for $i\ge 2$. By Lemma \ref{ellipsedistance}, since $f(b)<b$ for all $b\in (0,1/3)$, we have that $P_{b_i}A,Q_{b_i}B<\epsilon$ for all $1\le i\le m$. Define $V_i$ for $1\le i\le m$ to be 
\begin{align*}
      V_i=
  \bigg\{
  \begin{aligned}
B_{\delta}(Q_{b_i})\cap C, \text{ if } i \text{ is odd} \\
 B_{\delta}(P_{b_i})\cap C, \text{ if } i \text{ is even}
  \end{aligned}
\end{align*}
where $\delta$ is sufficiently small. We can certainly ensure that $V_i\subset B_{\epsilon}(A)\cup B_{\epsilon}(B)$ for all $1\le i\le m$ and that all the $V_i$ are disjoint. We can also ensure that $V_i$ is contained in the open upper right quadrant for all odd $i$ and $V_i$ is contained in the open lower left quadrant for all even $i$. Notice that by Lemma \ref{ellipsedistance} $f$ is continuous and strictly increasing on $(0,1/3)$ and hence has a continuous inverse $f^{-1}:f((0,1/3))\rightarrow (0,1/3)$. Let $U_m=V_m$. Recursively define 
\begin{align*}
      U_{m-i}=
  \bigg\{
  \begin{aligned}
  V_{m-i}\cap \{Q_{f(b)}\ |\ b\in (0,1/3), P_b\in U_{m-i+1}\}, \text{ if } m-i \text{ is odd} \\
 V_{m-i}\cap \{P_{f(b)}\ |\ b\in (0,1/3), Q_b\in U_{m-i+1}\}, \text{ if } m-i \text{ is even}
  \end{aligned}
\end{align*}
for $1\le i\le m-1$. By induction from $m$ to $1$ we have that for any $1\le i\le m$, $U_i$ is an open subset of $C$ such that if $i$ is even then $P_{b_i}\in U_i$, and if $i$ is odd then $Q_{b_i}\in U_i$. Our example that realizes $g$ will have that the points $X_j$ with $j\in A_i$ will lie in $U_i$. For any $b\in (0,1/3)$, define $R_b=(-2\sqrt{1-b^2},b),S_b=(2\sqrt{1-b^2},-b)$.\\

By Lemma \ref{ellipsedistance} and by symmetry we have that the unique farthest point on $C$ from $P_{b_{i+1}}$ is $Q_{b_i}$ and the unique farthest point on $C$ from $Q_{b_{i+1}}$ is $P_{b_i}$ for any $1\le i\le m-1$. This means that if $\delta$ is small enough, then for any $1\le i\le m-1$ and $1\le j\le m$ with $j\neq i$, we have that 
\begin{align*}
    EF>EG \text{ if }E\in V_{i+1},F\in V_i,G\in V_j\cup \{A,B\}.
\end{align*}
For all $j\in A_1$, pick $X_j$ to be an arbitrary point in $U_1$ so that the chosen points are distinct. We will now describe how to pick points $X_j\in U_i$ for $j\in A_i$ with $i\ge 2$ inductively on $i$. Suppose that for some $1\le i\le m-1$, we have chosen all the points $X_j\in U_i$ with $j\in A_i$ where the $X_j$ are distinct. Without loss of generality we may assume that $A_i=[l]$ for some positive integer $l$. We know that $j\in A_{i+1}$ if and only if $g(j)\in A_i$ so we can partition $A_{i+1}=g^{-1}(1)\cup g^{-1}(2)\cup \ldots \cup g^{-1}(l)$. From the definition of $U_i$, Lemma \ref{ellipsedistance} and by central symmetry of $C$ we know that there are points $Y_j\in U_{i+1}$ for $1\le j\le l$ such that the unique farthest point on $C$ from $Y_j$ is $X_j$. Now, let $\delta_i$ be sufficiently small and for each $j\in A_{i+1}$, pick $X_j$ to be an arbitrary point in $U_{i+1}\cap B_{\delta_i}(Y_j)$ so that all of the $X_j$ are disjoint. By picking $\delta_i$ to be small enough we can ensure that for each $j\in A_{i+1}$, the farthest point from $X_j$ among the points $\{X_{j'}\ |\ j'\in A_i\}$ is indeed $X_{g(j)}$. We can construct inductively $X_j$ for all $j\in S_A\backslash (A_1\cup \{a_0\})$. \\

Similarly, we can consider small neighbourhoods of $R_{b_i}$ for odd $i$ and $S_{b_i}$ for even $i$ to construct points $X_j$ for all $j\in S_B$ with similar properties. This would construct all the points $X_j$ for $j\in [n]$. By uniqueness of farthest points we can ensure that, by taking the neighbourhoods to be sufficiently small,T the farthest point from $X_j$ among $S=\{X_{j'}\ |\ j'\in [n]\}$ is indeed $X_{g(j)}$ for any $j\in S_A\cup S_B\backslash (A_1\cup B_1\cup \{a_0,b_0\})$. We now need to consider the case when $j\in A_1\cup B_1\cup \{a_0,b_0\}$.\\

Suppose that $j\in A_1$. By the triangle inequality  $|X_jX_{a_0}-Q_{b_1}A|\le \delta$. Note that $$Q_{b_1}A^2=b_1^2+4(1+\sqrt{1-b_1^2})^2=8-3b_1^2+8\sqrt{1-b_1^2}.$$
Let $b_0=f(b_1)$. By symmetry and Lemma \ref{ellipsedistance} we know that the farthest point on $C$ from $Q_{b_1}$ is $P_{b_0}$ and $b_0<b_1/3$, but we also know that as we move a point $E$ on $C$ from $P_{b_0}$ downwards (and right) the distance $Q_{b_1}E$ decreases at least until we get to the point $(0,-1)$. Thus, we know that $Q_{b_1}P_{b_i}<Q_{b_1}P_{b_2}$ for $i>2$. Notice that for any $i$, among points $Q_{b_i},P_{b_i},R_{b_i},S_{b_i}$ the (unique) farthest point from $Q_{b_1}$ is $P_{b_i}$ because it is farthest in both coordinates. Thus, based on our construction in order to show that for sufficiently small $\delta$, the farthest point from $X_j$ is indeed $X_{a_0}=A$ it is enough to show that $$Q_{b_1}P_{b_2},Q_{b_1}R_{b_1},Q_{b_1}B<Q_{b_1}A.$$
We have that \begin{align*}
Q_{b_1}B^2=b_1^2+4(1-\sqrt{1-b_1^2})^2, \\
Q_{b_1}A^2=b_1^2+4(1+\sqrt{1-b_1^2})^2, \\
Q_{b_1}R_{b_1}=4\sqrt{1-b_1^2}, \\
Q_{b_1}P_{b_2}^2=(b_1+b_2)^2+4(\sqrt{1-b_1^2}+\sqrt{1-b_2^2})^2.
\end{align*}
Clearly $Q_{b_1}B<Q_{b_1}A$. Since $1>\sqrt{1-b_1^2}$, we have that $Q_{b_1}A^2>4\cdot (2\sqrt{1-b_1^2})^2=Q_{b_1}R_{b_1}^2$ so $Q_{b_1}R_{b_1}<Q_{b_1}A$. From Lemma \ref{ellipsedistance} we also know that $b_1<b_2/3$ and hence we have that $$Q_{b_1}A^2-Q_{b_1}P_{b_2}^2=8\sqrt{1-b_1^2}-8\sqrt{1-b_1^2}\sqrt{1-b_2^2}+3b_2^2-2b_1b_2.$$
Since $3b_2^2>2b_1b_2$ and $\sqrt{1-b_2^2}<1$, we have that $Q_{b_1}P_{b_2}<Q_{b_1}A$. So if $\delta$ is sufficiently small, then for any $j\in A_1$, the farthest point from $X_j$ in $S$ is indeed $X_{a_0}=A$. Similarly, we can ensure that if $j\in B_1$, the unique farthest point from $X_j$ is $X_{b_0}=B$. Now, for any point $E=(x,y)\in C$ distinct from $B$, we have that 
$$AE^2=(x+2)^2+y^2=x^2+4x+4+1-\frac{x^2}{4}=\frac{3x^2+16x-44}{4}+16=\frac{(x-2)(3x+22)}{4}+16<16=AB^2.$$ So $X_{b_0}=B$ is the unique farthest point from $X_{a_0}=A$ among the points in $S$ and vice versa. Thus, we have that for all $j\in [n]$, the unique farthest point from $X_j$ in $S$ is $X_{g(j)}$. Now, similar to the argument at the end of the proof of Theorem \ref{Rkrealizable}, by slightly perturbing the points $X_1,X_2,\ldots X_n$ (not necessarily having them on $C$), we can also ensure that all the distances $X_iX_j$ are distinct. This completes the proof of the lemma.\QED

We are now ready to prove Theorem \ref{maxrealizable}.

\textit{Proof of Theorem \ref{maxrealizable}.} Suppose that $g$ is a max-realizable function. From the proof of Lemma \ref{no2cycles} we see that $G_g$ does not have cycles of length more than two. By Lemma \ref{partitiontree} there is a partition $G_g=F_0\cup F_1\cup \ldots \cup F_m$ for some $m\ge 0$ such that any edge in $G_g$ either goes from $F_i$ to $F_{i-1}$ for some $1\le i\le m$ or from $F_0$ to $F_0$. As in the proof of Lemma \ref{onecomponentrealizable}, $G_g[F_0]$ is a union of disjoint $2$-cycles. Let those cycles be $C_1,C_2,\ldots C_s$ for some $s\ge 1$. Define $H_i$ in the same way as in Lemma \ref{onecomponentrealizable}. Again as in the proof of that lemma we know that $H_i$ are disjoint and $H_g[H_i]$ are all disconnected from each other. Now consider the unit circle $S^1$. Let $A_1,A_2,\ldots A_s\in S^1$ and $B_1,B_2,\ldots B_s\in S^1$ be arbitrary such that all of the $2s$ points are distinct and $A_iB_i$ is a diameter of $S^1$ for all $i$. Now, let $\epsilon$ be sufficiently small. For any $1\le i\le s$, by Lemma \ref{onecomponentrealizable} we can select points $X_j\in B_{\epsilon}(A_i)\cup B_{\epsilon}(B_i)$ for all $j\in H_i$ such that the farthest point from $X_j$ among the points $\{X_{j'}\ |\ j'\in H_i\}$ is $X_{g(j)}$. 

Note that $A_iB_i>A_iB_j$ for any $i\neq j$. This means that we can certainly pick a sufficiently small $\epsilon$ so that for any $i\neq j$, if $E\in B_{\epsilon}(A_i),F\in B_{\epsilon}(B_i),G\in B_{\epsilon}(A_j)\cup B_{\epsilon}(B_j)$, then $EF>EG$. This will ensure that for any $j\in [n]$, the farthest point from $X_j$ among the points $\{X_{j'}\ |\ j'\in [n]\}$ is indeed $X_{g(j)}$. As before by a slight perturbation of the points $X_1,X_2,\ldots X_n$ we can also ensure that all the distances $X_iX_j$ are distinct. Thus, we have that $g$ is max-realizable in $\mathbb{R}^2$. This completes the proof. \QED
\\

If $(X,d)$ is a metric space, we will define a function $f:[n]\rightarrow [n]$ without fixed points to be \emph{min-realizable in }$(X,d)$ if there are points $A_1,A_2\ldots A_n\in X$ such that all distances $d(A_i,A_j)$ are distinct for all different pairs $\{i,j\}$ with $i\neq j$, and for any distinct $i,j\in [n]$, we have that $d(A_i,A_{f(i)})\le d(A_i,A_j)$. Define $f$ to be \emph{min-realizable} if it is min-realizable in some metric space. The above theorem tells us that $g$ is not very interesting and may suggest that a realizable pair of functions $(f,g)$ is realizable in $\mathbb{R}^2$ if and only if $f$ is min-realizable in $\mathbb{R}^2$. This is however not the case. In fact, one can check that every min-realizable function $f$ with $|\dom f|=6$ is also min-realizable in $\mathbb{R}^2$. However, there is a realizable pair $(f,g)$ that is not realizable in $\mathbb{R}^2$ with $|\dom f|=6$. We show this in the next proposition. We also note that, using the same method as in the proof of Theorem \ref{genrealizable}, we can show that the set of min-realizable functions is the same as the set of max-realizable functions and those are precisely the functions $f$ from $[n]$ to itself without fixed points such that $G_f$ has no cycles of length greater than two. 

\begin{prop}
    Let $f,g:[6]\rightarrow [6]$ be defined as follows: $f(i)=6$ for $i\in [5]$, $f(6)=1$, $g(i)=1$ for $i=2,3,4,5$, $g(1)=g(6)=2$. Then $(f,g)$ is a realizable pair that is not realizable in $\mathbb{R}^2$.
\end{prop}
\textit{Proof. } It is easy to see that $(f,g)$ is a nice pair and hence a realizable pair. Suppose that it is realizable in $\mathbb{R}^2$. Let $X_1,X_2,\ldots X_6\in \mathbb{R}^2$ be an example that realizes $(f,g)$, where $X_i$ corresponds to $i$ for all $i\in [6]$. We may assume that $X_6=O=(0,0)$ is the coordinate origin. We have that, by Proposition \ref{angle60}, for any distinct $i,j\in [5]$, $\angle X_iOX_j>\pi/3$. Suppose that ray $OX_i$ meets $S^1$ at $Y_i$ for $i\in [5]$, where $S^1$ is the unit circle. We may assume that $i_0=1$ and $i_0,i_1,\ldots i_4$ is a permutation of $[5]$ such that $Y_{i_0},Y_{i_1},\ldots Y_{i_4}$ appear on $S^1$ clockwise in that order. Let $T$ be the point on ray $OX_1$ such that $OT=OX_{i_4}$. Since $OX_1<OX_{i_4}$ we have $O-X_1-T$. Since $X_{i_1}O<X_{i_1}X_1$, we also know that $\angle X_{i_1}X_1O<\pi/2$ and hence $X_{i_1}X_1<X_{i_1}T$. Now, let $\alpha_j=\angle X_{i_j}OX_{i_{j+1}}$ for $j=0,1,\ldots 4$ where $i_5=i_0$. Since $\sum_{i=0}^4\alpha_i=2\pi$ and $\alpha_i>\pi/3$ for all $i$, we have that $\alpha_0+\alpha_4=2\pi-\alpha_1-\alpha_2-\alpha_3<\pi$. Now, by the cosine theorem, we have that $$X_{i_1}T^2=OX_{i_1}^2+OT^2-2OT\cdot OX_{i_1}\cos \alpha_0<OX_{i_1}^2+OX_{i_4}^2-2OX_{i_4}OX_{i_1}\cos (\alpha_0+\alpha_4)=X_{i_1}X_{i_4}^2.$$
But since $g(i_1)=1$, we have that $X_{i_1}X_1>X_{i_1}X_{i_4}>X_{i_1}T>X_{i_1}X_1$, which is a contradiction. This completes the proof.\QED
\\

Notice that if $k$ is fixed and $f$ is a function from $[n]$ to itself, then similar to the proof of the upper bound in Theorem \ref{Rkrealizable} if there is a vertex of $G_f$ with large enough in-degree (depending on $k$), we know that $f$ is not min-realizable in $\mathbb{R}^k$. However, having large in-degree is not the only obstacle. To show this, we will first define the \emph{complete upward binary tree of depth s} to be the tree $B^s$, where the vertex set is $\{(i,j)\ |\ 0\le i\le s,1\le j\le 2^i\}$ and the edge set consists of edges $(i+1,2j-1)\rightarrow (i,j),(i+1,2j)\rightarrow (i,j)$ for all $0\le i\le s-1,1\le j\le 2^i$. We define the \emph{level} of vertex $(i,j)\in B^s$ to be $i$. Notice that no vertex in $B^s$ has in-degree more than two. We now state the following theorem.
\begin{thm}\label{binarybounded}
    For any positive integer $k$, there is a large enough $s(k)$ such that the following holds. If $f$ is a function from $[n]$ to itself for some $n\in \mathbb{N}$, then for any $t\ge s(k)$ if $G_f$ has a copy of $B^t$. we have that $f$ is not min-realizable in $\mathbb{R}^{k}$.
\end{thm}
To prove this, we will first need a few lemmas. The following lemma gives a bound on the number of points in a ball such that no two points are close to each other. For any $A\in \mathbb{R}^k$ and positive real number $r$, let $B_r(A)$ denote the open ball in $\mathbb{R}^k$ with radius $r$ centered at $A$.
\begin{lem}\label{pointsinball}
    Let $k$ be a fixed positive integer. For any real $r>1$ and $A\in \mathbb{R}^k$, if $n>(2r+1)^k$ is a positive integer and $X_1,X_2,\ldots X_n\in B_r(A)$, then there are distinct $i,j$ such that $X_iX_j<1$.
\end{lem}
\textit{Proof.} Suppose that for some $r>1,A\in \mathbb{R}^k$, points $X_1,X_2,\ldots X_n\in B_r(A)$ are such that $X_iX_j\ge 1$ for all distinct $i,j\in [n]$. We have that $B_{1/2}(X_i)$ are all disjoint and are all subsets of $B_{r+1/2}(A)$. Let $\alpha$ denote the ($k$-dimensional) volume of an open ball in $\mathbb{R}^k$ with radius $1/2$. We have that $$(2r+1)^k\alpha=\Vol (B_{r+1/2}(A))\ge \Vol (\cup_{i=1}^nB_{1/2}(X_i))=n\alpha,$$
where the volumes are $k$-dimensional. So $n\le (2r+1)^k$, which completes the proof. \QED
\\

The next lemma will give a bound on the number of points that satisfy certain properties on minimal distances. It is essentially a generalization of the argument in the proof of the upper bound in Theorem \ref{Rkrealizable}.
\begin{lem}\label{starlikebound}
    Let $k$ be a positive integer. Then there is a real number $B_k\ge 1$ such that the following holds. If $P,X_1,X_2,\ldots X_n\in \mathbb{R}^k$ are such that $PX_i>6$ for all $i$ and $PX_i<X_iX_j+1$ for any distinct $i,j$, then $n\le B_k$.
\end{lem}
 For $k=1$ it is trivial so suppose that $k\ge 2$. Let $P,X_1,X_2,\ldots X_n\in \mathbb{R}^k$ be as above for some $n\in \mathbb{N}$. We may assume that $P$ is the coordinate origin. Now, let ray $PX_i$ meet $S^{k-1}$ at $Y_i$ for each $i$. We will show that for any distinct $i,j$, $\cos (\angle X_iPX_j)<2/3$. We may assume that $i=1,j=2$ and let $a=PX_1,b=PX_2,c=X_1X_2$. We may also assume that $a\ge b$. Note that $c>a-1$. Now we have
$$\cos (\angle X_iPX_j)=\frac{a^2+b^2-c^2}{2ab}<\frac{a^2+b^2-(a-1)^2}{2ab}=\frac{b^2+2a-1}{2ab}<\frac{1}{2}+\frac{1}{b}\le \frac{2}{3}$$
since $b\ge 6$. Thus, since $\angle Y_iPY_j=\angle X_iPX_j$, we have that $\angle Y_iPY_j>\arccos (2/3)$. Let $\alpha=\frac{\arccos (2/3)}{2}$. Let $v_i=\overrightarrow{PY_i}$ for each $i$. Now, by (\ref{7.5}) we have that the hyperspherical caps $C_{v_i,\alpha}$ in $S^{k-1}$ are all disjoint and hence 
$$A_{k-1}\ge \Area (\cup_{i=1}^n C_{v_i,\alpha})=\sum_{i=1}^n \Area (C_{v_i,\alpha})=nA_{k-1,\alpha}.$$
where the areas are $(k-1)$-dimensional. Thus, taking $B_k=\frac{A_{k-1}}{A_{k-1,\alpha}}$ gives the desired result. \QED
\\

Now we will prove Theorem \ref{binarybounded}.

\textit{Proof of Theorem \ref{binarybounded}}. Since $B^t\le B^{t+1}$, it is enough to show that there is a large enough $t$ such that if $f$ is a min-realizable function and $G_f$ has a copy of $B^t$, then $f$ is not min-realizable in $\mathbb{R}^k$. Assume the contrary. Suppose that $t\in \mathbb{N}$ is large. We can pick points $X_{i,j}\in \mathbb{R}^k$ with $0\le i\le t,j\in [2^i]$ such that all pairs of points have distinct distances and for any $0\le i\le t-1,j\in [2^i]$, $X_{i,j}$ is the closest point in $S=\{X_{i,j}\ |\  0\le i\le t,j\in [2^i]\}$ to $X_{i+1,2j-1}$ and the closest point in $S$ to $X_{i+1,2j}$. Suppose that $h:E(B^t)\rightarrow \mathbb{R}$ is defined such that if $e\in E(B^t)$ is an edge such that $e=(i_1,j_1)\rightarrow (i_2,j_2)$, then $h(e)=X_{i_1,j_1}X_{i_2,j_2}$. Also, define $g:B^t\backslash \{(0,1)\}\rightarrow B^t$ such that for any $0\le i\le t-1,1\le j\le 2^i$, $g((i+1,2j-1))=g((i+1,2j))=(i,j)$. Now we may sort the edges in $E(B^t)$ in order $e_1,e_2,\ldots e_{2^{t+1}-2}$ such that $h(e_1),h(e_2),\ldots h(e_{2^{t+1}-2})$ is a strictly increasing sequence. For convenience, let $p=2^{t+1}-2$. Now, for each $l\in [p]$, let $(i_l,j_l)$ be the head of edge $e_l$ (the edge starts from vertex $(i_l,j_l)$).\\

Suppose that $l$ is such that $i_l>1$. Then consider the directed path in $B^t$ $(i_l,j_l)=x\rightarrow g(x)\rightarrow g^{(2)}(x)\rightarrow \ldots g^{(i_l)}(x)=(0,1)$, where $g^{(s)}$ is $g$ applied $s$ times. By considering the point in $\mathbb{R}^k$ corresponding to $g^{(s)}(x)$ for $1\le s\le i_l-1$, we have that $h(g^{(s)}(x)g^{(s+1)}(x))<h(g^{(s-1)}(x)g^{(s)}(x))$. Therefore, in the sequence $e_1,e_2,\ldots e_p$ all of edges on the path from $(i_l,j_l)$ to $(0,1)$ other than $e_l$ appear before $e_l$. For convenience, let $P=X_{0,1}$. For any $1\le l\le p$, let $T_l$ be the graph consisting of edges $e_1,e_2,\ldots e_l$ and whose vertex set is the union of the vertices in $e_1,e_2,\ldots e_l$. Also, let $T_0$ be the graph with no edges and whose only vertex is $(0,1)$. From above we can see that $T_l$ has exactly one more edge and one more vertex than $T_{l-1}$ for $1\le l\le p$. We shall now prove the following claim.

\textbf{Claim.}  For any integers $1\le m<r\le p$, we have that if $r-m>(4t+1)^k$ then $h(e_r)\ge 2h(e_m)$.\\

\textit{Proof of Claim.} Suppose that $h(e_r)<2h(e_m)$. For each $1\le l\le r-m$, let $P_l=X_{i_{m+l},j_{m+l}}$. We have that the closest point to $P_l$ in $S$ is of distance $h(e_{m+l})$ to $P_l$ so for any $l,q$, we have that $h(e_m)<P_lP_q$. Suppose that $1\le l\le r-m$. Since $P_l$ corresponds to vertex $w=(i_{m+l},j_{m+l})$, consider $w\rightarrow g(w)\rightarrow \ldots g^{(i_{m+l})}(w)=(0,1)$. By above we know that $h(g^{(s)}(w)g^{(s+1)}(w))\le h(wg(w))=h(e_{m+l})\le h(e_r)<2h(e_m)$ for all $0\le s\le i_l-1$. By the triangle inequality
$$P_lP\le \sum_{s=0}^{i_l-1}h(g^{(s)}(w)g^{(s+1)}(w))<i_l2h(e_m)\le 2th(e_m).$$
Thus, we have that $P_l\in B_{2th(e_m)}(P)$ for all $l$. Applying Lemma \ref{pointsinball} scaled with a factor of $h(e_m)$ we have that $r-m\le (4t+1)^k$. This proves the claim. \QED

Consider the process of starting with an empty graph with vertex set the same as $B^t$ and adding edges $e_1,e_2,\ldots e_p$ one by one. We may assume that $t$ is odd and let $t=2d+1$. Suppose that $e_s$ is the first edge whose head is on level $d$. By above we know that $i_l<d$ for any $1\le l\le s-1$. Now suppose that $(i_s,j_s)=w_d\rightarrow w_{d-1}\rightarrow \ldots w_0=(0,1)$ is the directed path in $B^t$ from $(i_s,j_s)$ to $(0,1)$. For each $1\le l\le d$, let $y_l$ be the vertex in $B^t$ such that $y_l\rightarrow w_{l-1}$ is an edge in $B^t$ but $y_l\neq w_l$. Consider the sets 

$$S_l=\{v\in B^t\ | \text{ the directed path from }v\text{ to }(0,1) \text{ in }B^t\text{ passes through }y_l\}.$$
for $1\le l\le d$. The edges added so far have no vertex on a level greater than $d$. Thus, for each $l$, we can pick a vertex $v_l\in S_l$ of the highest level (that is still at most $d+1$) such that $v_l\rightarrow g(v_l)\not \in E(T_s)$ but $g(v_l)\in T_s$. For each $1\le l\le d$, let $H_l$ be the induced subgraph of $B^t$ on the set $\{g(v_l)\}\cup \{v\in B^t\ |\ g^{(i)}(v)=v_l \text{ for some }0\le i< d/B_k-1\}$, where $B_k$ is defined as in Lemma \ref{starlikebound}. By construction $H_l$ are edge-disjoint graphs. From the claim we see that adding at least $(4t+1)^k+1$ edges means that we have to at least double the $h$ function. Applying this $t$ times means that adding at least $(4t+1)^{k+1}>t((4t+1)^k+1)$ edges we have to multiply our $h$ function by at least $2^t$ times. Let $s_i=s+i(4t+1)^{k+1}$ for $0\le i \le t$. Notice that $s<2^{d+1}$ and hence $t(4t+1)^{k+1}+s<2^t+2^{d+1}<p$ for sufficiently large $t$ since exponential beats polynomial. This means that $s_t<p$ for large $t$.\\

Define $H_{l,i}=H_l\cap T_{s_i}$ for $0\le i<d/B_k,l\in [d]$. Notice that $H_{l,0}$ consists of just the vertex $g(v_l)$. Now, say that $H_l$ is $i$-good for $1\le i<d/B_k-1,l\in [d]$ if for any vertex $v\in H_l\backslash H_{l,i-1}$ such that $g(v)\in H_{l,i-1}$, we have that $v\in H_{l,i}$. Essentially this means that all the neighbouring vertices to $H_{l,i-1}$ are added to $H_{l,i}$. Now suppose that $i$ is fixed. We will show that there are at most $B_k$ of the graphs $H_l$ that are not $i$-good. Let $S\subset [d]$ be the set of indices such that for each $l\in S$, we have that $H_l$ is not $i$-good. That means that for each $l\in S$, there is some vertex $x_l\in H_l\backslash H_{l,i-1}$ such that $g(x_l)\in H_{l,i-1}$ and $x_l\not \in H_{l,i}$. Now suppose that $Q_l,R_l$ are the points in $\mathbb{R}^k$ corresponding to $x_l,g(x_l)$ respectively. By definition of $s_i$ we have that $h(e_m)>2^th(e_{s_{i-1}})$ for any $m\ge s_i$ and therefore $Q_lR_l>2^th(e_{s_{i-1}})$. Since $g(w_l)\in H_{l,i-1}$, we know that similar to above, by the triangle inequality, $R_lP\le th(e_{s_{i-1}})$. For large enough $t$, we have that $2^t>7t$ and for any disjoint $l_1,l_2\in S$, we have that $$Q_{l_1}Q_{l_2}>Q_{l_1}R_{l_1}\ge Q_{l_1}P-PR_{l_1}\ge Q_{l_1}P-th(e_{s_{i-1}})$$ and we also have that for any $l\in S$ $$Q_lP\ge Q_lR_l-R_lP>6th(e_{s_{i-1}}).$$ Thus, we can apply Lemma \ref{starlikebound} scaled by a factor of $th(e_{s_{i-1}})$ to conclude that $|S|\le B_k$.\\

Let $m$ be the largest integer such that $m<d/B_k<t$. Clearly $m\ge d/B_k-1$ and by above we have that for any $i$, at most $B_k$ of the $H_l$ are not $i$-good. Therefore, there must be some $l'\in [d]$ such that $H_{l'}$ is $i$-good for all $1\le i\le m$. Notice that, by induction on $j$, if $H_l$ is $i$-good for $1\le i\le j$, then we must have that $H_{l,j}$ must contain all vertices in the set $\{g(v_l)\}\cup \{v\in H_l\ |\ g^{(j')}(v)=v_l \text{ for some }0\le j'\le j-1\}$. This means that $H_{l'}=H_{l',m}$. From $e_s$ to $e_{s_m}$ we have added at most $m(4t+1)^{k+1}\le d(8d+5)^{k+1}$ edges. But we have added all the edges of $H_{l'}$, which there are more than $2^{m-1}>2^{d/B_k-2}>d(8d+5)^{k+1}$ for large enough $d$. This gives a contradiction which completes the proof of the theorem. \QED
\\

We have therefore shown that having a large in-degree is not the only obstacle to being min-realizable in $\mathbb{R}^k$. Ideally, we want to find, if possible, a stronger condition for being min-realizable in $\mathbb{R}^k$ or at least in $\mathbb{R}^2$. Therefore, we leave with the following two open problems:
\begin{prob}
    For any given $k$, classify all the functions that are min-realizable in $\mathbb{R}^{k}$.
\end{prob}
\begin{prob}
    For any given $k$, classify all the realizable pairs of functions that are realizable in $\mathbb{R}^k$.
\end{prob}

\bibliographystyle{plain}
\bibliography{arxiv}

\begin{thebibliography}{1}

\bibitem{10.5555/7228}
B.~Bollob\'{a}s.
\newblock {\em Combinatorics: set systems, hypergraphs, families of vectors, and combinatorial probability}.
\newblock Cambridge University Press, 1986.

\bibitem{Croft}
H.~T. Croft.
\newblock personal communication, 2022.

\bibitem{hyperspherical}
S.~Li.
\newblock Concise formulas for the area and volume of a hyperspherical cap.
\newblock {\em Asian Journal of Mathematics \& Statistics}, 4, 01 2011.

\end{thebibliography}

\end{document}